\documentclass[10pt,twoside]{article}

\usepackage{graphicx}
\usepackage{fancyhdr}
\renewcommand{\Pi}{\pi}

\usepackage [usenames] {color}
\usepackage [colorlinks=true,
linkcolor=webgreen,
filecolor=webbrown,
citecolor=webgreen] {hyperref}
\definecolor {webgreen} {rgb} {0,.5,0}
\definecolor {webbrown} {rgb} {.6,0,0}
                        
\usepackage {color}
\usepackage {float}

\usepackage{amsfonts}
\usepackage{amstext}
\usepackage{amssymb,amsmath,amsbsy}

\usepackage{hyperref}

\newtheorem{theorem}{Theorem}

\newtheorem{remark}{Remark}
\newtheorem{definition}{Definition}
\newtheorem{lemma}{Lemma}

\newcommand{\proof}[2]{\noindent {\bf Proof} {\bf #1}\,\,\,{#2}}

\newcounter{thcount}
\def\qued{\hfill \vrule height 8pt width6pt depth0pt}

\def\R{\Bbb R}

\def\0{\bold 0}
\def\1{\bold 1}
\def\cc{\bold c}
\def\dd{\bold d}
\def\a{\bold a}
\def\b{\bold b}
\def\q{\bold q}
\def\r{\bold r}
\def\x{\bold x}
\def\uu{\bold u}
\def\vv{\bold v}

\def\z{\bold z}

\def\y{\bold y}

\def\eps{\varepsilon}
\def\X{\bold X}

\def\V{\bold V}

\def\A{\bold A}
\def\Y{\bold Y}
\def\I{\bold I}

\def\L{\bold L}
\def\M{\bold M}

\def\RR{\bold R}
\def\W{\bold W}
\def\BB{\bold B}

\def\DD{\bold D}
\def\E{\bold E}
\def\tr{\text{tr}\,}
\def\diag{\text{diag}\,}

\def\dist{\text{dist}}

\def\Span{\text{Span}\,}

\def\Vol{{\text{Vol}}\,}
\def\part{\cal P} 
\def\tu{\tilde \uu}

\begin{document}


\begin{center}
{\Large {\bf Beyond the Expanders}}
\end{center}

\begin{center}
{\large {\bf Marianna~Bolla}}\\
\end{center}

\begin{center}
{Institute of Mathematics, Budapest University of Technology and Economics, \\
Egry Jozsef u. 1, 1111 Budapest, Hungary}; E-mail: {\tt {marib@math.bme.hu}}\\
\end{center}


\date{}

\renewcommand\abstractname{Abstract}
                                             
\begin{abstract}
\noindent \\
Expander graphs are widely used in communication problems and construction
of error correcting codes.
In such graphs, information gets through very
quickly. Typically, it is not true for social or biological networks,
though we may find a partition of the vertices such that the induced subgraphs
on them and the bipartite subgraphs between any pair of them 
exhibit regular  behavior of information flow within or between
the vertex subsets. 
Implications between spectral and regularity properties are discussed.


\end{abstract}

\section {Introduction}
\label{intro}

We want to go beyond the expander graphs that -- for four decades --
 have played  an important role 
in communication networks; for a summary, see e.g., Chung~\cite{Chung0} and
Hoory et al.~\cite{Hoory}.
Roughly speaking, the expansion property means that each subset of the graph's
vertices has ``many'' neighbors (combinatorial view), and hence, information 
gets through such a
graph very ``quickly'' (probabilistic view). 
We will not give  exact definitions of expanders
here as those contain many parameters which are not used later.
We rather refer to the spectral and random walk characterization of 
such graphs,
as discussed, among others by Alon~\cite{Alon0}, and  
Meila and Shi~\cite{Meila}. 

The general framework of an edge-weighted graph will be used.
Expanders have a spectral gap bounded away from zero, where -- for a connected
graph -- this gap is 
defined as the minimum distance between the 
normalized Laplacian spectrum (apart from the trivial zero eigenvalue) and the
endpoints of the [0,2] interval, the possible range of the spectrum.
The larger the spectral gap, the more our graph resembles a random graph and
exhibits quasi-random properties, e.g.,
the edge densities within any subset and between any two subsets of its 
vertices do not differ too much of what is expected, see the Expander Mixing 
Lemma~\ref{lem} of Section~\ref{pre}. Quasi-random properties and spectral
gap of random graphs with given expected degrees are discussed 
in Chung and Graham~\cite{Chung}, and Coja-Oghlan and Lanka~\cite{Coja}.

However, the spectral gap appears not at the ends of the normalized Laplacian
spectrum in case of 
generalized random or 
generalized quasi-random graphs that, in the presence of
 $k\ge 2$ underlying clusters, have $k$ 
eigenvalues (including the 
zero) separated from 1, while the bulk of the spectrum is located around 1,
see e.g.,~\cite{Bolla3}. These structures are usual in social or biological
networks having $k$ clusters of vertices (that belong to social groups or
similarly functioning enzymes) such that the edge density within the clusters
and between any pair of the clusters is homogeneous. 

Our conjecture is that $k$ so-called structural eigenvalues (separated from 1)
in the normalized Laplacian spectrum are indications of
such a structure, while the near 1 eigenvalues are responsible for the
pairwise regularities. The clusters themselves can be recovered by applying
the $k$-means algorithm for the vertex representatives obtained by the
eigenvectors corresponding to the structural eigenvalues (apart from the zero). 
For the $k=2$ case
we will give an exact relation between the eigenvalue separation (of the 
non-trivial structural eigenvalue from the bulk of the spectrum) and the 
volume regularity
of the cluster pair that is obtained by the $k$-means algorithm applied for
the coordinates of the transformed eigenvector belonging to the non-trivial
structural eigenvalue, see Theorem~\ref{th} of Section~\ref{tetel}. 
To eliminate the trivial eigenvalue-eigenvector pair,
we shall rather use the normalized modularity spectrum of~\cite{Bolla5} 
 that plays an important
role in finding the extrema of some penalized versions of the Newman-Girvan
modularity introduced in~\cite{New1}.
Theorem~\ref{thk} of Section~\ref{anova} gives an estimation for the 
extent of volume-regularity of the different cluster pairs in the $k>2$
case based on the spectral gap and the $k$-variance of the vertex
representatives.

In~\cite{Coja1,Mc}, the authors give algorithms -- based on low rank
approximation -- to find a regular partition if $k$ is known and our
graph comes from a generalized random graph model with $k$ clusters.
Without knowing $k$, there are constructions -- like~\cite{Frieze} --
based on refinement of partitions and leading to a very fine partition
with number of clusters depending merely on the constant
 ruling the regularity of the cluster pairs. On the contrary, our purpose
is to estimate the extent of the regularity of the cluster pairs by means
of spectral gaps and  eigenvectors. The
estimations given are relevant only in the presence of a large spectral gap
(between some structural and the other eigenvalues) and special 
classification properties of the
eigenvectors corresponding to the structural eigenvalues, 
see Theorem~\ref{thk} of Section~\ref{anova}. 
In this case,
the algorithm is straightforward via $k$-means clustering.

\section{Preliminaries and statement of purpose}
\label{pre}

Let $G=(V, \W)$
be a graph on $n$ vertices, where the $n\times n$ symmetric matrix $\W$ has
non-negative real entries and zero diagonal. 
  Here $w_{ij}$ is the similarity 
between vertices $i$ and $j$, where 0 similarity means no connection/edge
at all. A simple graph is a special case of it with 0-1 weights.
Without loss of generality 
\begin{equation}\label{feltetel}
   \sum_{i=1}^n \sum_{j=1}^n w_{ij}=1
\end{equation}
will be supposed.
Hence, $\W$ is a joint distribution, with marginal entries
$$
 d_i =\sum_{j=1}^n w_{ij} , \quad i=1,\dots ,n
$$
which are the \textit{generalized vertex degrees}
collected  in the main diagonal of the diagonal \textit{degree matrix} 
$\DD =\diag (\dd )$, \, $\dd =(d_1 ,\dots ,d_n )^T$. 
In~\cite{Bolla,Bolla1}  we investigated the spectral gap of the
\textit{normalized Laplacian}
 $\L_D =\I -\DD^{-1/2} \W \DD^{-1/2}$,
where $\I$ denotes the identity matrix of appropriate size.

Suppose that our graph is connected ($\W$ is irreducible).
Let $0=\lambda_1 <\lambda_2 \le \dots \le \lambda_n \le 2$  denote
the eigenvalues of the symmetric normalized Laplacian $\L_D$ with 
corresponding
unit-norm, pairwise orthogonal eigenvectors $\uu_1 ,\dots ,\uu_n$.
Namely, $\uu_1 =(\sqrt{d_1},\dots , \sqrt{d_n})^T =\sqrt{\dd }$.
In the random walk setup $\DD^{-1} \W$ is the transition matrix (its entry
in the $(i,j)$-th position is the conditional probability of moving from
vertex $i$ to vertex $j$ in one step, given that we are in $i$)  which
is a stochastic matrix with eigenvalues $1-\lambda_i$
and corresponding eigenvectors $\DD^{-1/2} \uu_i$  $(i=1,\dots ,n )$.
``Good'' expanders have a $\lambda_2$ bounded away from zero, that also implies
the separation of the isoperimetric number 
\begin{equation}\label{h}
  h(G) =\min_{U\subset V : \,\Vol (U) \le \frac12} 
         \frac{w (U,{\overline U})}{\Vol (U)} ,
\end{equation}
where for $X,Y\subset V$: $w(X,Y) =\sum_{i\in X} \sum_{j\in Y} w_{ij}$ is the
weighted cut between $X$ and $Y$,  while
$\Vol (U) =\sum_{i\in U} d_i$ is the volume of $U\subset V$. In view 
of~(\ref{feltetel}), $\Vol (V) =1$, this is why the minimum is taken on 
vertex sets having volume at
most $\frac12$. In~\cite{Bolla1}, we proved that
\begin{equation}\label{iso}
 \frac12 \lambda_2 \le h(G) \le \min \{ 1 ,\sqrt{2\lambda_2} \} ,
\end{equation}
while in the $\lambda_2 \le 1$ case the stronger upper estimation 
$$
 h (G) \le \sqrt{\lambda_2 (2-\lambda_2 )}
$$
holds. (We remark that $\lambda_2 \le \frac{n}{n-1}$ always holds.)

If a network does not have a ``large'' $\lambda_2$ (compared to the natural
lower bound), or equivalently -- in view of the above  inequalities -- it has
a relatively ``small'' isoperimetric number, then the 2-partition 
of the vertices giving the minimum in~(\ref{h}) 
indicates a bottleneck, or equivalently, a low conductivity edge-set between
two disjoint vertex clusters such that the random walk gets through with small 
probability between them,
but -- as some equivalent notions will indicate -- it is rapidly mixing within
the clusters. To find the clusters, the coordinates of the transformed 
eigenvector $\DD^{-1/2} \uu_2$ will be used. In~\cite{Bolla}, we proved that
for the weighted 2-variance of this vector's coordinates 
\begin{equation}\label{2szoras}
 {S }_2^2 (\DD^{-1/2} \uu_2 ) \le \frac{\lambda_2 }{\lambda_3} 
\end{equation}
holds. For a general $2\le k\le n$, the notion of $k$-variance
-- in the Analysis of Variance sense -- is the following.
The weighted $k$-variance of the $k$-dimensional vertex representatives 
$\x_1 ,\dots , \x_n$ comprising the row vectors of the $n\times k$ matrix $\X$ 
is defined by
\begin{equation}\label{kszoras}
 {S}_k^2 (\X ) =\min_{P_k \in {\part}_k }
  {S}_k^2 (P_k ,\X ) =\min_{P_k =(V_1 ,\dots ,V_k )}
\sum_{a=1}^k \sum_{j\in V_a } d_j \| \x_j -{ \cc }_a \|^2 ,
\end{equation}
where ${\cc }_a =\frac1{\Vol (V_a ) } \sum_{j\in V_a } d_j \x_j $ is the
weighted center of cluster $V_a$  $(a=1,\dots ,k )$ and ${\part}_k$ denotes
the set of $k$-partitions of the vertices.
We remark that 
${S }_2^2 (\DD^{-1/2} \uu_1 , \DD^{-1/2} \uu_2 )={S }_2^2 (\DD^{-1/2} \uu_2 )$,
since
 $\DD^{-1/2} \uu_1 =\1$ is the all 1's vector.

The above results were generalized for minimizing the normalized $k$-way cut
\begin{equation}\label{szam}
 f_k (P_k ,G) = \sum_{a=1}^{k-1} \sum_{b=a+1}^k \left( \frac1{\Vol (V_a)}+
   \frac1{\Vol (V_b)} \right) w (V_a ,V_b )=
  k-  \sum_{a=1}^{k} \frac{w (V_a ,V_a )}{\Vol (V_a )} 
\end{equation}
of the $k$-partition $P_k =(V_1 ,\dots ,V_k )$ over the set of all possible
$k$-partitions. Let
$$
 f_k (G) =\min_{P_k \in {\part}_k } f_k (P_k ,G) 
$$
be the \textit{minimum normalized k-way cut} of the underlying weighted 
graph $G =(V, \W )$. In fact, $f_2 (G)$ is the symmetric version of the
isoperimetric number and $f_2 (G)\le 2h(G)$. In ~\cite{Bolla1} 
we proved that
\begin{equation}\label{fk}
  \sum_{i=1}^k \lambda_i \le f_k (G) \le c^2 \sum_{i=1}^k \lambda_i ,
\end{equation}
where the upper estimation is relevant only in the case when 
$S_k^2 (\uu_1 ,\dots ,\uu_k )$ is small enough and the constant $c$
depends on this minimum $k$-variance of the vertex representatives.

The \textit{normalized Newman-Girvan modularity} is defined in~\cite{Bolla5}
as the penalized version of the Newman-Girvan modularity~\cite{New1}
in the following way.  The normalized $k$-way modularity of $P_k =
(V_1 ,\dots ,V_k )$ is
\begin{equation}\label{normod}
\begin{aligned}
 Q_k (P_k ,G) &=  
 \sum_{a=1}^k \frac1{\Vol (V_a) }  \sum_{i,j\in V_a} (w_{ij} -d_i d_j ) =
   \sum_{a=1}^k \frac1{\Vol (V_a) } [ w(V_a ,V_a) - \Vol^2 (V_a )] \\
  &= \sum_{a=1}^k \frac{w(V_a ,V_a)}{\Vol (V_a )} -1 =k-1-f_k (P_k )  ,
\end{aligned}
\end{equation}
and
$$
 Q_k  (G)=\max_{P_k \in {\part}_k } Q_k (P_k ,G) 
$$
is the \textit{maximum normalized k-way Newman-Girvan modularity} 
of the underlying weighted graph $G =(V, \W )$. For given $k$, maximizing 
this modularity is equivalent to minimizing the normalized cut and can be solved
by the same spectral technique. 
In fact, it is more convenient to use the spectral decomposition
of the normalized modularity matrix $\BB_D = \I -\L_D -\sqrt{\dd} \sqrt{\dd}^T$
with eigenvalues $\beta_1 \ge \dots \ge \beta_n$, that are the
numbers $1-\lambda_i$ 
with eigenvectors $\uu_i$  $(i=2,\dots ,n)$
and the zero with corresponding unit-norm eigenvector $\sqrt{\dd}$.
In~\cite{Bolla1,Bolla5}, we also show that a spectral gap between $\lambda_k$
and $\lambda_{k+1}$ is an indication of $k$ clusters with low inter-cluster
connections; further, the intra-cluster connections ($w_{ij}$) between 
vertices $i$ and $j$
of the same cluster are higher 
than expected under the hypothesis of independence 
(in view of which the vertices are connected
with probability $d_i d_j $). In the random walk framework, the random walk
stays within the clusters with high probability.

Conversely, minimizing the above modularity will result in clusters with
high inter- and low intra-cluster connections. 
In~\cite{Bolla5}, we proved that 
\begin{equation}\label{Qmin}
 \min_{P_k \in {\part}_k } Q_k (P_k ,G) \ge \sum_{i=1}^k \beta_{n+1-i} . 
\end{equation}
The existence of $k$ ``large''
(significantly larger than 1) eigenvalues in the normalized Laplacian 
spectrum, or equivalently, the existence of $k$ negative eigenvalues
(separated from 0) in the normalized modularity spectrum is an  indication of 
$k$ clusters with the above 
property. In the random walk setup: the walk stays within the clusters with
low probability. 

These two types of network structures are frequently called
community or anti-community structure. 
These are the two extreme cases, when $f_k (P_k ,G)$ is either minimized or
maximized, and the optimization gives $k$ clusters with  either strong
intra-cluster and weak inter-cluster connections, or vice versa.
Some networks exhibit a more general,
still regular behavior: the  vertices can be classified into  $k$ clusters
such that the information-flow within them and between any pair of them is 
homogeneous. In terms of random walks, the walk stays within clusters or
switches between clusters with probabilities characteristic for the cluster
pair. That is, if the random walk moves from a vertex of cluster $V_a$ 
to a vertex
of cluster $V_b$, then the probability of doing this does not depend on the 
actual vertices,
it merely depends on their cluster memberships, $a,b=1,\dots ,k$.

In this context, we examined the following  generalized random graph model,
that corresponds to the ideal case:
given the number of clusters 
$k$, the vertices of the graph independently belong to the clusters;
 further, conditioned on the cluster
memberships, vertices $i\in V_a$ and $j\in V_b$ are connected with
probability $p_{ab}$, independently of each other, $1\le a,b\le k$.
Applying the results~\cite{Bolla3} for the spectral characterization 
of some
noisy random graphs, we are able to prove that the normalized modularity 
spectrum of   a generalized random graph is the following:
 there exists a positive  number $\theta <1$, independent
of $n$, such that 
there are exactly $k-1$ so-called structural eigenvalues of $\BB_D$ that are
greater than $\theta -o(1)$,
while all the others are o(1) in absolute value. It is equivalent that
$\L_D$ has $k$ eigenvalues (including the zero) separated from 1.

The $k=1$ case corresponds to quasi-random graphs and the above characterization
corresponds to the eigenvalue separation of such graphs, 
discussed in~\cite{Chung}. The authors
also prove some implications between the so-called quasi-random properties.
For example, for dense graphs, ``good'' eigenvalue separation is equivalent
to ``low'' discrepancy (of the induced subgraphs' densities from the overall 
edge density). 

For the $k\ge 2$ case, generalized quasi-random graphs were introduced by
Lov\'asz and T. S\'os~\cite{LovSos}. These graphs are 
deterministic counterparts
of generalized random graphs with the same spectral properties. In fact, the
authors define so-called generalized quasi-random graph sequences by means
of graph convergence that also implies the convergence of spectra.
Though, the spectrum itself does not carry enough information for the
cluster structure of the graph, together with some classification
properties of the structural eigenvectors it does.
We want to prove some implication between the spectral gap and
the volume-regularity of the cluster pairs, also using the structural
eigenvectors.

The notion of volume regularity was introduced by Alon et al.~\cite{Alon}.
We shall use a slightly modified version of this notion.

\begin{definition} 
Let $G=(V, \W)$  be weighted graph with $\Vol (V) =1$.
The disjoint pair $(A,B)$ is
$\alpha$-\textit{volume regular}
 if for all $X\subset A$, $Y\subset B$ 
we have
\begin{equation}\label{jeles}
| w (X, Y) -\rho (A,B) \Vol (X) \Vol (Y)| \le \alpha \sqrt{\Vol (A) \Vol (B)} ,
\end{equation}
where $\rho (A,B) =\frac{e(A,B)}{ \Vol (A) \Vol (B)}$ is the relative
inter-cluster density of $(A,B)$. 
\end{definition}

Our definition was inspired by the Expander Mixing Lemma stated e.g., 
in~\cite{Hoory} for regular graphs and in~\cite{Chung0} for simple
graphs in the context of quasi-random properties.  
Now we formulate it for  edge-weighted graphs 
on a general degree sequence. We also include the proof as a preparation
for the proof of Theorem~\ref{th} of Section~\ref{tetel}. 

\begin{lemma}\label{lem} 
(\textbf{Expander Mixing Lemma for Weighted Graphs})
Let $G=(V,W)$ be a weighted graph and suppose that
$\Vol (V)=1$.  Then for all $X,Y\subset V$:
$$
\begin{aligned}
 | w (X, Y) - \Vol (X) \Vol (Y)| &\le \| \BB_D \|\cdot 
 \sqrt{\Vol(X) (1-\Vol (X)) \Vol(Y) (1-\Vol (Y))}  \\
 &\le \| \BB_D \|\cdot  \sqrt{\Vol (X) \Vol (Y)} ,
\end{aligned}
$$
where $\| \BB_D \|$ is the spectral norm of the normalized modularity
matrix of $G$.
\end{lemma}

\proof {}{
Let $X\subset A$, $Y\subset B$ 
and $\1_U \in \R^n$ denote the indicator vector of $U\subset V$. Further,
$\x :=\DD^{1/2} \1_X$ and $\y :=\DD^{1/2} \1_Y$.

We use the spectral decomposition $\DD^{-1/2} \W \DD^{-1/2} =\sum_{i=1}^n \rho_i
\uu_i \uu_i^T$ , where $\rho_i =1-\lambda_i$ $(i=2,\dots ,n)$ are 
eigenvalues of $\BB_D$ and $\rho_1 =1$ with corresponding unit-norm 
eigenvector $\uu_1 =\sqrt{\dd }= \DD^{1/2} \1$. We remark  that $\uu_1$ is also
an eigenvector of $\BB_D$ corresponding to the eigenvalue zero, hence
$\| \BB_D \| =\max_{i\ge 2} |\rho_i |$. 
Let $\x = \sum_{i=1}^n x_i \uu_i$ and
$\y = \sum_{i=1}^n y_i \uu_i$ be the expansions of $\x$ and $\y$ in the 
orthonormal basis $\uu_1 ,\dots ,\uu_n$ with coordinates
$x_i =\x^T \uu_i$ and $y_i =\y^T \uu_i$, respectively. Observe that
 $x_1 =\Vol (X)$, $y_1 =\Vol (Y)$ and
$\sum_{i=1}^n x_i^2 =\| \x \|^2 =\Vol (X)$,
$\sum_{i=1}^n y_i^2 =\| \y \|^2 =\Vol (Y)$.
Based on these,
$$
\begin{aligned}
 | w (X, Y) - \Vol (X) \Vol (Y)| &=|\sum_{i=2}^n \rho_i x_i y_i | \le
 \| \BB_D \| \cdot |\sum_{i=2}^n x_i y_i |  \\
 &\le \| \BB_D \| \cdot \sqrt{\sum_{i=2}^n x^2_i \sum_{i=2}^n y^2_i } \\
 &\le \| \BB_D \| \cdot \sqrt{\Vol(X) (1-\Vol (X)) \Vol(Y) (1-\Vol (Y))}  \\
 &\le  \| \BB_D \| \cdot \sqrt{\Vol(X) \Vol(Y)} ,
\end{aligned}
$$
where we also used the triangle and the Cauchy-Schwarz inequalities.
\qued\\ }

We remark that the spectral gap of $G$ is  $1-\| \BB_D \|$, hence 
-- in view of Lemma~\ref{lem} -- the density between any 
two subsets of ``good'' expanders is near to what is expected. 
On the contrary, in the above definition
of volume regularity, the $X,Y$ pairs are disjoint, 
and a  ``small'' $\alpha$ indicates that the $(A,B)$ pair is like
a bipartite expander, see e.g.,~\cite{Chung0}.  

In the next section we shall prove the following statement for the $k=2$ case:
if one eigenvalue jumps out of 
the bulk of the normalized modularity spectrum, 
then clustering the coordinates
of the corresponding transformed eigenvector into 2 parts 
(by minimizing the 2-variance of its coordinates)
will result in an $\alpha$-volume regular partition of the vertices, where
$\alpha$ depends on the spectral gap.

We may go further: if $k-1$ (so-called structural) eigenvalues jump out
of the normalized modularity spectrum, then clustering the
representatives of the vertices -- obtained by the corresponding eigenvectors
in the usual way -- into $k$ clusters will result in $\alpha$-volume 
regular pairs, 
where $\alpha$ depends on the spectral gap (between the structural
eigenvalues and the bulk of the spectrum) and the $k$-variance
of the vertex representatives based on the eigenvectors corresponding to the
structural eigenvalues. In Section~\ref{anova}, we give an estimation for
$\alpha$ in the $k\ge 2$ case; further, we extend the estimation to  the
clusters themselves.

\section{Eigenvalue separation and volume regularity (k=2 case)}
\label{tetel}

\begin{theorem}\label{th}
Let $G=(V, \W)$ is an edge-weighted graph on $n$ vertices, with generalized
degrees $d_1 ,\dots ,d_n$ and $\DD =\diag (d_1 ,\dots ,d_n )$. Suppose that
$\Vol (V) =1$. 
Let the eigenvalues  of $\DD^{-1/2} \W \DD^{-1/2}$,
enumerated in decreasing absolute values, be 
$$
 1=\rho_1 >|\rho_2 |=\theta >\eps \ge |\rho _i | , \quad i\ge 3 .
$$
The partition $(A,B)$ of $V$ is defined so that it minimizes
the weighted 2-variance of the coordinates of $\DD^{-1/2} \uu_2$, 
where $\uu_2$ is the unit-norm eigenvector
belonging to $\rho_2$. Then the $(A,B)$ pair is 
${\cal O} (\sqrt{\frac{1-\theta }{1-\eps}})$-volume regular.
\end{theorem}

\proof {}{
We  use the notations of Lemma~\ref{lem}'s proof.  
Let $X\subset A$, $Y\subset B$. Fort short,
$\x :=\DD^{1/2} \1_X$, $\y :=\DD^{1/2} \1_Y$,
 $\a :=\DD^{1/2} \1_A$, $\b :=\DD^{1/2} \1_B$.
With $\rho :=\rho (A,B)$ and $\M := \W -\rho \dd \dd^T$,
\begin{equation}\label{egy}
 | w (X, Y) -\rho \Vol (X) \Vol (Y)| = |\1_X^T \M \1_Y| = |\x^T 
 (\DD^{-1/2} \W \DD^{-1/2} -\rho \sqrt{\dd} \sqrt{\dd}^T) \y | .
\end{equation}
Using the spectral decomposition $\DD^{-1/2} \W \DD^{-1/2} =\sum_{i=1}^n \rho_i
\uu_i \uu_i^T$ and the fact that $\uu_1 =\sqrt{\dd }= \DD^{1/2} \1$, we can 
write~(\ref{egy}) as
\begin{equation}\label{ketto}
|(1-\rho )x_1 y_1 + \rho_2 x_2 y_2 + \sum_{i=3}^n \rho_i x_i y_i | ,
\end{equation}
where $\x = \sum_{i=1}^n x_i \uu_i$ and
$\y = \sum_{i=1}^n y_i \uu_i$ is the expansion of $\x$ and $\y$ in the 
orthonormal basis $\uu_1 ,\dots ,\uu_n$ with coordinates
$x_i =\x^T \uu_i$ and $y_i =\y^T \uu_i$, respectively.

First we will prove that $1-\rho$ is governed by $\rho_2$;
more precisely,  $|1-\rho |\le |\rho_2 |+\eps$.
Applying the arguments of Lemma~\ref{lem} and the above formulas for  the 
special  $A,B\subset V$ yields
\begin{equation}\label{sign}
\begin{aligned}
 \Vol (A) \Vol (B)\cdot ( \rho  - 1 ) &= w(A,B) - \Vol (A) \Vol (B) =\\
 &=\a^T (\DD^{-1/2} \W \DD^{-1/2} -\sqrt{\dd} \sqrt{\dd}^T ) \b  
 =\rho_2 a_2 b_2 + \sum_{i=3}^n \rho_i a_i b_i ,
\end{aligned}
\end{equation}
where  $\a = \sum_{i=1}^n a_i \uu_i$ and
$\b = \sum_{i=1}^n b_i \uu_i$ is the expansion of $\a$ and $\b$ in the 
orthonormal basis $\uu_1 ,\dots ,\uu_n$, respectively. The separation of $A$
and $B$ is based on the vector $\DD^{-1/2} \uu_2$ which  has both 
negative and positive coordinates, since  $\uu_2$
is orthogonal to $\uu_1$ of all positive coordinates. 
With formulas, $\a +\b =\uu_1$, and hence, $a_2 +b_2 = \uu_1^T \uu_2 =0$.
(If it is the eigenvalue $\lambda_2$ of the normalized Laplacian that 
is the farthest 
from 1, then the corresponding eigenvector, our $\uu_2$, is also 
called ``Fiedler-vector'' as the  two-partition of the vertices 
into  two loosely  connected parts was based on
the signs of its coordinates in the early paper of Fiedler~\cite{Fiedler}). 
If $\theta$ is
much larger than $\eps$, the first term in the last formula
of~(\ref{sign}) 
-- apart from a term of ${\cal O} (|\eps |)$ -- will dominate the sign of 
$\rho  - 1$ which is therefore opposite to the sign of $\rho_2$. 

Therefore, we will distinguish between two cases.
\begin{itemize}
\item
If $\lambda_2 <1-\eps$, then $\rho_2 =1-\lambda_2 > \eps >0$, and
in view of the inequalities between
the minimum normalized cut and the smallest positive normalized 
Laplacian eigenvalue (apply~(\ref{fk}) for the $k=2$ case):
\begin{equation}\label{iso}
 \rho \ge f_2 (G)=\min_{U\subset V} \frac{w (U, {\bar U} )}{\Vol (U) 
  \Vol({\bar U })} \ge \lambda_2 =1-\rho_2  ,
\end{equation}
therefore $1-\rho \le \rho_2$, as $1-\rho$ is also positive due to the
considerations before. Further, the estimation, due to~(\ref{2szoras}),
\begin{equation}\label{egyik}
{S}_2^2 (\DD^{-1/2} \uu_2 ) \le \frac{\lambda_2}{\lambda_3} 
 \le \frac{1-\rho_2}{1-\eps} = \frac{1-\theta }{1-\eps }
\end{equation}
also follows.

\item 
If $1-\eps \le \lambda_2 \le \frac{n}{n-1}$,
then -- provided $\frac1{n-1}\le \eps$ -- 
it is the eigenvalue $\lambda_n$ that is the farthest from 1, and hence,
greater than $1+\eps$. Consequently, 
$-\eps < \rho_2 =1-\lambda_n <0$, and hence,  by~(\ref{normod}) and
(\ref{Qmin}):
$$
 \rho_2 +\rho_{neg} \le Q_2 ((A,B), \W ) =(2-1)-f_2 ((A,B), \W) =1-\rho ,
$$ 
where $\rho_{neg} =\min \{1-\lambda_{n-1} , 0 \}$, and $|\rho_{neg} | < \eps$.
Note, that in this case $1-\rho$ is negative that yields 
$|1-\rho| \le |\rho_2 | +\eps$.
Now the optimum $A,B$ is obtained by minimizing the 2-variance of the
coordinates of the transformed eigenvector $\DD^{-1/2} \uu_2$ (now $\uu_2$
belongs to $\lambda_n$ and $\rho_2$ at the same time) for which
the following relation -- like (\ref{2szoras}) -- can be proved: 
\begin{equation}\label{un}
{S}_2^2 (\DD^{-1/2} \uu_2 ) = {\cal O}(\frac{2-\lambda_n}{2-\lambda_{n-1}})
 ={\cal O} (\frac{\beta_n +1}{\beta_{n-1}+1} ) 
  = {\cal O} (\frac{1-\theta }{1-\eps } ) ,
\end{equation}
where $\beta$'s are eigenvalues of the normalized modularity matrix.
Indeed, in lack of dominant vertices, there is a relation 
between the largest and smallest normalized Laplacian eigenvalues of $G$
and  $\overline G$, respectively,
where the complement graph ${\overline G} =(V, {\overline \W})$ is defined
such that ${\overline w}_{ij} =1-w_{ij}$ $(i\ne j)$ and ${\overline w}_{ii}=0$
$(i=1,\dots ,n )$.

If the two largest absolute value eigenvalues of the normalized modularity
matrix are of different sign, then we are able to find a gap at least
$\theta -\eps$ between eigenvalues of the same sign.
\end{itemize}



 Therefore,  (\ref{ketto}) can be estimated from above with
\begin{equation}\label{error}
 |\rho_2 | \cdot |x_1 y_1 + x_2 y_2 | + \eps x_1 y_1 +
\max_{i\ge 3} |\rho_i |\cdot |\sum_{i=3}^n x_i y_i |  .
\end{equation}
As for the second term, $\eps x_1 y_1 =\eps \Vol (X) \Vol (Y)$, 
so it does not need further treatment.

Using the Cauchy-Schwarz inequality, the last term can be estimated from above
with
$$
 \eps \sqrt{\sum_{i=3}^n x^2_i \sum_{i=3}^n y^2_i } 
  \eps \sqrt{\sum_{i=2}^n x^2_i \sum_{i=2}^n y^2_i } 
\le
 \eps \sqrt{\Vol(X) (1-\Vol (X)) \Vol(Y) (1-\Vol (Y))} 
 \le \eps \sqrt{\Vol(X) \Vol(Y)}  ,
$$ 
since $x_1 =\Vol (X)$, $y_1 =\Vol (Y)$ and
$\sum_{i=1}^n x_i^2 =\| \x \|^2 =\Vol (X)$,
$\sum_{i=1}^n y_i^2 =\| \y \|^2 =\Vol (Y)$.
 
The first term is reminiscent of an equation for the coordinates of
orthogonal vectors. Therefore, we project the vectors $\uu_1$, $\uu_2$ onto
the subspace $F=\Span \{ \a ,\b \}$. In fact, $\uu_1 =\a +\b$, and hence,
$\uu_1 \in F$. The vector $\uu_2$ can be decomposed as
\begin{equation}\label{u2}
 \uu_2 = \frac{\uu_2^T \a }{\Vol (A)} \a + \frac{\uu_2^T \b }{\Vol (B)} \b
      +\q ,
\end{equation}
where $\q$ is the component orthogonal to $F$.
For the squared distance $\| \q \|^2$ between $\uu_2$ and $F$, in~\cite{Bolla},
we proved that it is equal to the weighted 2-variance 
${S}_2^2 (\DD^{-1/2} \uu_2 )$ 
and in~(\ref{egyik})  
we estimated it from above with $\frac{1-\theta }{1-\eps }$.
(In the $\rho_2 =1-\lambda_n$ case similar upper estimation 
works using~(\ref{un})). Let $s^2$ denote
this minimum $2$-variance of the coordinates of $\DD^{-1/2} \uu_2$
 (in both cases).

To estimate
$
 a_1 b_1 +a_2 b_2 =(\uu_1^T \a )(\uu_1^T \b )+(\uu_2^T \a )(\uu_2^T \b ) 
$,
the problem is that the pairwise orthogonal vectors $\uu_1 ,
\uu_2$ and $\a ,\b$ are not in the same subspace of $\R^n$ as, in general,
$\uu_2 \notin F$. However, by an argument proved in~\cite{Bolla}, we can
find orthogonal, unit-norm vectors ${\tu}_1 ,{\tu}_2 \in F$ such that
\begin{equation}\label{vari}
 \| \uu_1 -\tu_1 \|^2 +  \| \uu_2 -\tu_2 \|^2 \le 2s^2 ,
\end{equation}
where, in view of $\uu_1 \in F$, $\tu_1 =\uu_1$. Let $\r :=\uu_2 -\tu_2 $. Since
$\tu_1^T \a ,\tu_2^T \a$ and $\tu_1^T \b ,\tu_2^T \b$ are coordinates
of the orthogonal vectors $\a ,\b$ in the basis $\tu_1 ,\tu_2$,
$$
 (\tu_1^T \a )(\tu_1^T \b )+(\tu_2^T \a )(\tu_2^T \b )  =0 ,
$$
and because of $\tu_2^T \a +\tu_2^T \b = \tu_2^T \uu_1 =0$, 
$$
 \tu_2^T \a =-\tu_2^T \b =\sqrt{\Vol (A) \Vol (B )} =:c .
$$
Therefore,
$$
\begin{aligned}
 &|(\uu_1^T \a )(\uu_1^T \b )+(\uu_2^T \a )(\uu_2^T \b ) | 
 =|\Vol (A) \Vol (B )+[(\tu_2 +\r )^T \a ][(\tu_2 +\r )^T \b ] | \\
 &=|\Vol (A) \Vol (B )+[c +\r^T \a ][-c +\r^T \b ] | 
= |c (-\r^T \a +\r^T \b )+(\r^T \a)(\r^T \b ) | \\
&\le |c| \sqrt{ \| \r \|^2 \| \b -\a \|^2 } + 
  \sqrt{\| \r \|^2 \| \a \|^2 }\sqrt{\| \r \|^2 \| \b \|^2 }  \\   
 & \le\sqrt{\Vol (A) \Vol (B )} (|\r \| +\| \r \|^2 ) \le
\sqrt{\Vol (A) \Vol (B )} (\sqrt2 s+2s^2 ) ,
\end{aligned}
$$
using  (\ref{vari}) and the fact that $\| \b -\a \|^2 =1$.

Now we estimate
$x_1 y_1 + x_2 y_2 = (\uu_1^T \x )(\uu_1^T \y )+(\uu_2^T \x )(\uu_2^T \y )$.
Going back to~(\ref{u2}) we have
$$
 \uu_2^T \x = \frac{\uu_2^T \a}{\Vol (A)} \a^T \x +
              \frac{\uu_2^T \b}{\Vol (B)} \b^T \x + \q^T \x =
              \frac{\Vol (X)}{\Vol (A)} \uu_2^T \a + \q^T \x  ,
$$
and similarly,
$$
 \uu_2^T \y = \frac{\uu_2^T \a}{\Vol (A)} \a^T \y +
              \frac{\uu_2^T \b}{\Vol (B)} \b^T \y + \q^T \y =
              \frac{\Vol (Y)}{\Vol (B)} \uu_2^T \b + \q^T \y ,
$$
that in view of $\| \q \|^2 =s^2$ yields
$$
\begin{aligned}
 &x_1 y_1 +x_2 y_2 =|(\uu_1^T \x )(\uu_1^T \y )+(\uu_2^T \x )(\uu_2^T \y )|=\\
 &|\Vol (X) \Vol (Y) + (\frac{\Vol (X)}{\Vol (A)} \uu_2^T \a +\q^T \x ) 
                     (\frac{\Vol (Y)}{\Vol (B)} \uu_2^T \b +\q^T \y ) |  \\
&\le |\Vol (X) \Vol (Y) + (\frac{\Vol (X)}{\Vol (A)} \uu_2^T \a )  
                          (\frac{\Vol (Y)}{\Vol (B)} \uu_2^T \b )| \\
&+|(\q^T \x ) (\frac{\Vol (Y)}{\Vol (B)} \uu_2^T \b ) +
  (\q^T \y ) (\frac{\Vol (X)}{\Vol (A)} \uu_2^T \a ) + (\q^T \x )(\q^T \y )| \\
&\le \frac{\Vol (X)}{\Vol (A)}\frac{\Vol (Y)}{\Vol (B)} 
 |\Vol (A) \Vol (B) + (\uu_2^T \a ) ( \uu_2^T \b ) | \\
&+ \| \q \| \| \x \| \frac{\Vol (Y)}{\Vol (B)} \| \uu_2 \| \| \b \| 
 + \| \q \| \| \y \| \frac{\Vol (X)}{\Vol (A)} \| \uu_2 \| \| \a \| 
 + \| \q \|^2 \| \x \| \| \y \| \\
&\le \sqrt{\Vol (A) \Vol (B)} (\sqrt2 s+2s^2 ) +
 \| \q \| \sqrt{\Vol (X)} \frac{\Vol (Y)}{\Vol (B)} \sqrt{\Vol (B)} \\
&+\| \q \| \sqrt{\Vol (Y)} \frac{\Vol (X)}{\Vol (A)} \sqrt{\Vol (A)} + 
 \| \q \|^2 \sqrt{\Vol (X)} \sqrt{\Vol (Y)} \\
&=\sqrt{\Vol (A) \Vol (B)} (\sqrt2 s +2s^2 ) +
 \| \q \| \sqrt{\Vol (X)} \sqrt{\Vol (Y)} 
          (\frac{\sqrt{\Vol (Y)}}{\sqrt{\Vol (B)}} +
           \frac{\sqrt{\Vol (X)}}{\sqrt{\Vol (A)}} + \| \q \| ) \\
&\le \sqrt{\Vol (A) \Vol (B)}  [(\sqrt2 s+2 s^2 +s (2+s)] =
     \sqrt{\Vol (A) \Vol (B)}  [(\sqrt2 +2 )s +3s^2 ]  \\
& \le
  \sqrt{\Vol (A) \Vol (B)}  (\sqrt2 +5 )s   .
\end{aligned}
$$
Summarizing, the second and third terms in~(\ref{error}) are estimated from
above with  $\eps \sqrt{\Vol (X) \Vol (Y)} \le\eps \sqrt{\Vol (A) \Vol (B)}$.
Because of $\eps <\theta$, by an easy calculation it
follows that it is less than $\sqrt{\frac{1-\theta }{1-\eps}}$.
Therefore, the constant $\alpha$ of the $(A,B)$ pair's regularity is 
${\cal O} (\sqrt{\frac{1-\theta}{1-\eps}})$.
\qued\\ }

\begin{remark}
{
The statement has relevance only if $\theta$ is much larger than $\eps$. In 
this case the spectral gap between the largest absolute value eigenvalue and
the others in the normalized modularity  spectrum indicates a regular
2-partition of the graph that can be constructed based on the eigenvector
belonging to the structural eigenvalue.
}
\end{remark}



\section{Analysis of Variance setup (the $k > 2$ case)}\label{anova}

\begin{theorem}\label{thk}
Let $G=(V, \W)$ is an edge-weighted graph on $n$ vertices, with generalized
degrees $d_1 ,\dots ,d_n$ and $\DD =\diag (d_1 ,\dots ,d_n )$. Suppose that
$\Vol (V) =1$. 
Let the eigenvalues  of $\DD^{-1/2} \W \DD^{-1/2}$,
enumerated in decreasing absolute values, be 
$$
 1=\rho_1 >|\rho_2 | \ge \dots \ge |\rho_k | >\eps \ge |\rho _i | , 
 \quad i\ge k+1 .
$$
The partition $(V_1,\dots ,V_k )$ of $V$ is defined so that it minimizes
the weighted k-variance of the vertex representatives
obtained as row vectors of the $n\times k$ matrix $\X$ of column vectors
 $\DD^{-1/2} \uu_i$, 
where $\uu_i$ is the unit-norm eigenvector
belonging to $\rho_i$  $(i=1,\dots ,k)$. 
With the notation $s^2 =S_k^2 (\X )$, the $(V_i, V_j)$ pairs are
$2(\sqrt{2}s +\eps)$-volume regular $(i\ne j)$ and for the
clusters $V_i$ $(i=1,\dots ,k )$  the following holds:
for all $X,Y\subset V_i$ 
we have that
\begin{equation}\label{jelese}
| w (X, Y) -\rho (V_i ) \Vol (X) \Vol (Y)| \le 2(\sqrt{2}s +\eps) \Vol (V_i ) ,
\end{equation}
where $\rho (V_i ) =\frac{w(V_i ,V_i )}{ \Vol^2 (V_i ) }$ is the relative
intra-cluster density of $V_i$. 
\end{theorem}

\proof {}{
Denoting by $\uu_1 ,\dots ,\uu_k$ the eigenvectors belonging to the so-called
structural eigenvalues $\rho_1 ,\dots ,\rho_k$, the
representatives $\r_1 ,\dots ,\r_n$ of the vertices are row vectors   
of the matrix $\X =(\x_1 ,\dots, \x_k )$, where $\x_i =\DD^{-1/2} \uu_i$
$(i=1,\dots ,k)$ and
the trivial $\x_1 =\1$ (belonging to $\rho_1 =1$) can be omitted, 
see~(\ref{kszoras}).
The  minimum $k$-variance $S_k^2 (\X )$ 
of the $k$-dimensional (actually, $(k-1)$-dimensional) 
representatives  is as small as $s^2$. 
Suppose that the minimum   $k$-variance is attained by the
$k$-partition $(V_1 ,\dots ,V_k )$ of the vertices. 

By an easy analysis of variance argument of~\cite{Bolla1,Bolla3} it follows that
$$
 s^2 =\sum_{i=1}^k \dist^2 (\uu_i , F ) ,
$$
where $F =\Span \{ \DD^{1/2} \z_1 , \dots ,\DD^{1/2} \z_k \}$ with the so-called
normalized partition vectors $\z_1 ,\dots ,\z_k$ of coordinates
$z_{ji} = \frac1{\Vol (V_i )}$ if $j\in V_i$ and 0, otherwise $(i=1,\dots ,k)$. 
Note that the vectors $\DD^{1/2} \z_1 , \dots ,\DD^{1/2} \z_k$ form an 
orthonormal system. By~\cite{Bolla,Bolla1} we can find another orthonormal
system $\vv_1 ,\dots ,\vv_k \in F$ such that
$$
  \sum_{i=1}^k \| \uu_i -\vv_i \|^2 \le 2s^2 .
$$
With these vectors, we construct the following $k$-rank approximation of
the matrix $\DD^{-1/2} \W \DD^{-1/2} =\sum_{i=1}^n \rho_i \uu_i \uu_i^T$: it is
approximated by $\sum_{i=1}^k \rho_i \vv_i \vv_i^T$ with the following
accuracy (in spectral norm):
\begin{equation}\label{Frobenius}
 \| \sum_{i=1}^n \rho_i \uu_i \uu_i^T -\sum_{i=1}^k \rho_i \vv_i \vv_i^T \|
 \le \sum_{i=1}^k |\rho_i | \cdot \| \uu_i \uu_i^T -\vv_i \vv_i^T \| +
 \| \sum_{i=k+1}^n \rho_i \uu_i \uu_i^T \| \le 
 \sqrt{\sum_{i=1}^k \sin^2 \sigma_i }+\eps \le \sqrt{2}s +\eps ,
\end{equation}
where $\sigma_i$ is the 
angle between $\uu_i$  and $\vv_i$, and for it,
$\sin \frac{\sigma_i}2 = \frac12 \|u_i -v_i \|$ holds, therefore
$$
 \sin^2 \sigma_i =(2 \sin \frac{\sigma_i}2 \cos \frac{\sigma_i}2 )^2 =
 \frac14 \| u_i -v_i \|^2 (4 - \| u_i -v_i \|^2 )  , \quad i=1,\dots ,k .
$$
Hence, the above difference can be estimated from above with
$\sqrt{2}s +\eps$ in spectral norm.

Based on these considerations and the fact that the cut norm is less than or
equal to  the spectral norm,
the densities to be estimated in the defining formula~(\ref{jeles}) of
volume regularity can be written in terms of stepwise constant vectors in
the following way. The vectors $\y_i := \DD^{-1/2} \vv_i$ are stepwise 
constants on the partition $(V_1 ,\dots ,V_k )$, $i=1,\dots ,k$. The matrix 
$\sum_{i=1}^k \rho_i \y_i \y_i^T $ is therefore a symmetric block-matrix on $k\times k$
blocks belonging to the above partition of the vertices. Let ${\tilde w}_{ab}$
denote its entries in the $(a,b)$ block $(a,b=1,\dots ,k)$. 
Using~(\ref{Frobenius}), the following approximation of the matrix $\W$ is
performed:
$$
\| \W - \DD (\sum_{i=1}^k \rho_i \y_i \y_i^T ) \DD \| =
 \| \DD^{1/2} ( \DD^{-1/2} \W \DD^{-1/2} -\sum_{i=1}^k \rho_i \vv_i \vv_i^T ) 
 \DD^{1/2} \|
 \le \| \DD \|^{1/2} (\sqrt{2}s +\eps ) \| \DD \|^{1/2} .
$$
Therefore, the entries of $\W$ -- for $i\in V_a$, $j\in V_b$ -- 
can be decomposed as
$$
  w_{ij} = d_i d_j {\tilde w}_{ab} +\eta_{ij} ,
$$
where the cut norm and spectral norm of 
the $n\times n$ symmetric error matrix $\E =(\eta_{ij} )$   is at most
$\| \DD \| (\sqrt{2}s +\eps)$. 
But we will restrict the error matrix to $V_a \times V_b$:
its entries are $\eta_{ij}$'s for $i\in V_a , j\in V_b$, and zeros otherwise.
Denoting the restricted matrix by $\E^{ab}$, and the restricted diagonal
matrices by $\DD^a$ and $\DD^b$, respectively, the following finer estimation
holds:
$$
 \| \E^{ab} \|_{\square} \le |\ \DD^a \|^{1/2} \cdot \| \DD^b \|^{1/2} \cdot 
  (\sqrt{2}s +\eps ) \le
 \sqrt{\Vol (V_a ) \Vol (V_b )} (\sqrt{2}s +\eps ) .
$$
Consequently, for $a,b=1,\dots ,k$:
$$
\begin{aligned}
 &| w (X, Y) -\rho (V_a ,V_b ) \Vol (X) \Vol (Y)|  = 
\left| \sum_{i\in X} \sum_{j \in Y} (d_i d_j {\tilde w}_{ab}+\eta^{ab}_{ij}) - 
\frac{\Vol (X) \Vol (Y )}{\Vol (V_a)  \Vol (V_b )}
\sum_{i\in V_a } \sum_{j \in V_b } (d_i d_j {\tilde w}_{ab} +\eta^{ab}_{ij} )
   \right|  \\
&=\left| \sum_{i\in X} \sum_{j \in Y} \eta^{ab}_{ij} - 
\frac{\Vol (X) \Vol (Y )}{\Vol (V_a)  \Vol (V_b )}
\sum_{i\in V_a } \sum_{j \in V_b } \eta^{ab}_{ij} \right| \le 
 2 (\sqrt{2}s +\eps )  \sqrt{\Vol (V_a ) \Vol (V_b )} ,
\end{aligned}
$$
that gives the required statement both in the $a\ne b$ and $a=b$ case.



}

\begin{remark}
{
In the $k=2$ case, the estimate of Theorem~\ref{th} has the same order of 
magnitude as that of Theorem~\ref{thk}, since 
$s^2 ={\cal O} (\sqrt{\frac{1-\theta }{1-\eps}})$.
The statement has only relevance for an integer $k \in [2,n)$ such that 
there is a remarkable spectral gap
between $\theta :=|\rho_k |$ and $|\rho_{k+1} |$ in the normalized modularity
spectrum, i.e., the  so-called structural eigenvalues $\rho_1 , \dots ,\rho_k$
 are far apart from
zero, while the others are in an $\eps$ distance from zero, in absolute value.
This is a necessary condition for $s^2$ to be ``small''.
 As it is not
sufficient, instead of $\theta$ and $\eps$, 
the estimation of Theorem~\ref{thk} is given in terms of $s$ and $\eps$.
Indeed, by perturbation results of
spectral subspaces for symmetric matrices~\cite{Bath}, 
$s^2$ itself can be estimated from above by the spectral gap 
between the
$k$ structural and the other eigenvalues when $\rho_2 ,\dots ,\rho_k$ have
the same sign (the situation of strong community or anti-community structure). 
}
\end{remark}

\section*{Acknowledgement}

The author wishes to thank Vera T. S\'os, L\'aszl\'o Lov\'asz, and
Mikl\'os Simonovits for their useful advices.

\end{document}